\documentclass[11pt]{amsart}
\usepackage{verbatim, graphicx, rotating}
\usepackage{mdwlist}
\usepackage{enumerate,mdwlist}
\usepackage{url}
\usepackage{amsthm}
\usepackage{hyperref}

\usepackage {amssymb}

\input xy
\xyoption{all}
\input epsf
\newlength{\tabwidth}
\newlength{\tabheight}
\newlength{\tabrule}
\newlength{\tabwidthx}
\newlength{\tabheightx}

\def\gentabbox#1#2#3#4{\vbox to \tabheight{\setlength{\tabrule}{#3}%
  \setlength{\tabwidthx}{#1\tabwidth}\addtolength{\tabwidthx}{\tabrule}%

\setlength{\tabheightx}{#2\tabheight}\addtolength{\tabheightx}{-\tabheight}%
  \hbox to #1\tabwidth{%
    \hspace{-0.5\tabrule}\rule{\tabrule}{#2\tabheight}\hspace{-\tabrule}%
    \vbox to #2\tabheight{\hsize=\tabwidthx%
      \vspace{-0.5\tabrule}\hrule width\tabwidthx height\tabrule%
      \vspace{-0.5\tabrule}\vfil%
      \hbox to \tabwidthx{\hss#4\hss}%
        \vfil\vspace{-0.5\tabrule}%
      \hrule width\tabwidthx height\tabrule\vspace{-0.5\tabrule}}%
    \hspace{-\tabrule}\rule{\tabrule}{#2\tabheight}\hspace{-0.5\tabrule}}%
  \vspace{-\tabheightx}}}
\def\genblankbox#1#2{\vbox to \tabheight{\vfil\hbox to
#1\tabwidth{\hfil}}}
\def\tabbox#1#2#3{\gentabbox{#1}{#2}{0.4pt}{\strut #3}}

\catcode`\:=13 \catcode`\.=13 \catcode`\;=13 
\catcode`\>=13 \catcode`\^=13
\def:#1\\{\hbox{$#1$}}
\def.#1{\tabbox{1}{1}{$#1$}}
\def>#1{\tabbox{2}{1}{$#1$}}
\def^#1{\tabbox{1}{2}{$#1$}}
\def;{\genblankbox{1}{1}\relax}
\catcode`\:=12 \catcode`\.=12 \catcode`\;=12 
\catcode`\>=12 \catcode`\^=7

\newcommand{\field}{\mathbb}
\newcommand{\liealgebra}{\mathfrak}
\newcommand{\la}{\liealgebra}

\newcommand{\C}{{\field C}}

\newcommand{\Z}{{\field Z}}
\newcommand{\N}{{\field N}}

\renewcommand{\b}{\liealgebra b}

\newcommand{\n}{{\la n}}

\newcommand{\ga}{\alpha}

\newtheorem{prop}{Proposition}

\newtheorem{lemma}[prop]{Lemma}

\newtheorem{conjecture}[prop]{Conjecture}

\theoremstyle{definition}

\newtheorem{remark}[prop]{Remark}
\newtheorem{example}[prop]{Example}

\newtheorem{definition}[prop]{Definition}

\begin{document}
\title[Schubert Calculus Conjectures]
{Some Conjectures Regarding certain Schubert Structure Constants in Lie Types $B$ and $D$}

\author{Benjamin J. Wyser}
\date{\today}

%\thanks{}
%\address{}
%\email{}

\maketitle

In this note, I detail some conjectures related to Schubert calculus in types $B$ and $D$.  These are mentioned in \cite{Wyser-12b}, but are not spelled out in a precise way there due to the fairly large amount of notation required to give precise statements.  They are based upon an as yet conjectural understanding of the weak order poset of $L$-orbits on $G/B$, where $(G,L) = (SO(2n+1,\C),\C^* \times SO(2n-1,\C))$ (type $B$) or $(SO(2n,\C),\C^* \times SO(2n-2,\C))$ (type $D$).  In each case, $L$ is a Levi subgroup of $L$ which is ``spherical", meaning it acts with finitely many orbits on the flag variety $G/B$.  The idea is to identify the $L$-stable Richardson varieties with certain of the conjectural orbit parameters.

Fix the following notation:
\begin{itemize}
	\item $G$ a complex reductive algebraic group;
	\item $B,B^- \subseteq G$ opposite Borel subgroups;
	\item $T = B \cap B^-$, a maximal torus of $G$;
	\item $W = N_G(T)/T$, the Weyl group;
	\item $P,P^- \subseteq G$ opposite parabolic subgroups containing $B,B^-$, respectively;
	\item $L = P \cap P^-$ the common Levi factor of $P,P^-$.
\end{itemize}
	
For each $w \in W$, there exists a \textit{Schubert class} $S_w = [\overline{B^-wB/B}] \in H^*(G/B)$.  It is well-known that the classes $\{S_w\}_{w \in W}$ form a $\Z$-basis for $H^*(G/B)$.  As such, for any $u,v \in W$, we have
\[ S_u \cdot S_v = \displaystyle\sum_{w \in W} c_{u,v}^w S_w \]
in $H^*(G/B)$, for uniquely determined integers $c_{u,v}^w$.  These integers are the \textit{Schubert structure constants}.

It remains an open problem, even in type $A$, to give a positive combinatorial formula for an arbitrary Schubert constant $c_{u,v}^w$.  Various special cases are understood; see \cite{Wyser-12b} for one example of such a special case, as well as the other works referenced in the introduction to that paper.

The rule of \cite{Wyser-12b} is deduced using the following simple observation:  In the notation above, if $L$ is a spherical Levi subgroup, $X_u = \overline{BuB/B}$ is a Schubert variety stable under $P$, and $X^v = \overline{B^-vB/B}$ is an opposite Schubert variety stable under $P^-$, then the \textit{Richardson variety} $X_u^v := X_u \cap X^v$ is stable under $L$.  Since $L$ has finitely many orbits on $G/B$, it of course has finitely many orbits on $X_u^v$, and so there is a dense $L$-orbit on $X_u^v$.  Thus the Richardson variety $X_u^v$ is the closure of an $L$-orbit.  Since knowing the Schubert constants $c_{w_0u,v}^w$ is the same as knowing the expansion in the Schubert basis of $[X_u^v]$, one can use a theorem of M. Brion to deduce these constants.  This theorem expresses the class of any spherical subgroup orbit closure on $G/B$ as a sum of Schubert cycles, in terms of weighted paths in the weak order graph for the set of $L$-orbits on $G/B$.

To turn these simple observations into a useful combinatorial rule, one must do the following:
\begin{enumerate}
	\item Determine exactly which Richardson varieties are stable under $L$.
	\item Have a concrete combinatorial model of the set $L \backslash G/B$, as well as an understanding of its weak order in terms of that model.
	\item Match up the Richardson varieties stable under $L$ with the appropriate $L$-orbit closures.
\end{enumerate}

Step (1) above is easy --- the $L$-stable Richardson varieties are of the form $X_u^v$ where $u$ (resp. $v$) is a maximal (resp. minimal) length coset representative of $W_P \backslash W$, with $W_P$ the parabolic subgroup of $W$ associated to $P$.  In \cite{Wyser-12b}, the remaining steps are carried out for the pair $(G,L) = (GL(p+q,\C),GL(p,\C) \times GL(q,\C))$ using the known parametrization of $L$-orbits in that case, those results being due to Matsuki-Oshima and Yamamoto (\cite{Matsuki-Oshima-90,Yamamoto-97}).

For the type $B$ and $D$ pairs mentioned above, no parametrization of the $L$-orbits has appeared in the literature.  What \textit{has} appeared (again, cf. \cite{Matsuki-Oshima-90}) is a parametrization of $K$-orbit closures, where $K$ is the symmetric subgroup $S(O(2,\C) \times O(2n-1,\C))$ (type $B$) or $S(O(2,\C) \times O(2n-2,\C))$ (type $D$).  Note that each of these symmetric subgroups is disconnected, having two components, and that $L=K^0$ in both cases.  Thus each $K$-orbit is either disconnected, being a union of two distinct $L$-orbits, or is connected and coincides with a single $L$-orbit.  A parametrization of the $L$-orbits then amounts to knowing which $K$-orbits are connected and which are not.  Given this information, one can then easily deduce the weak order on $L$-orbits from the (known) weak order on $K$-orbits.

Below, I give specific conjectures on these matters.  The conjectures on the orbit sets were suggested first to me by Peter Trapa, so I attribute them to him.  They have been verified using ATLAS (available at \url{http://www.liegroups.org/}) through reasonably high rank.  Assuming those conjectures to be valid, we carry out the remaining steps for the type $B$ and $D$ pairs mentioned above, identifying the $L$-stable Richardson varieties, matching them with the proper combinatorial invariants for the $L$-orbit closures, and getting a rule for structure constants from this matching.  This portion of the argument can be proven just as we do for the type $A$ case in \cite{Wyser-12b}, although I do not write down the proofs here.  What makes them conjectures at this point is simply the lack of a proof of Trapa's conjectures regarding the set of $L$-orbits.

The Schubert calculus conjectures themselves have been tested using C code written by Dominic Searles which implements the algorithm of \cite{Knutson-03}.  Conjecture \ref{conj:type-b-schubert} (type $B$) has been verified to hold through rank $6$, and in rank $7$ for all suitable $u,v$ with $l(u) + l(v) \leq 12$.  Conjecture \ref{conj:type-d}  (type $D$) has been verified to hold through rank $7$.

\section*{Details}
As indicated above, we discuss the type $B$ pair $(SO(2n+1,\C),\C^* \times SO(2n-1,\C))$ and the type $D$ pair $(SO(2n,\C),\C^* \times SO(2n-2,\C))$.  In each case, $L$ is the Levi factor of the maximal parabolic subgroup $P$ of $G$ corresponding to omission of the simple root $\ga_1 = x_1 - x_2$, and is known to be spherical.  We realize $W$ as signed permutations of $\{1,\hdots,n\}$ (changing an even number of signs in type $D$), and denote these in one-line notation with bars over some of the numbers to indicate negative values.  For instance, $2 \overline{1} \overline{3}$ denotes the permutation sending $1$ to $2$, $2$ to $-1$, and $3$ to $-3$.

In type $B$, the following facts are easy to check:
\begin{enumerate}
	\item Left cosets $W_P \backslash W$ consist of all signed permutations whose one-line notations have either $1$ or $\overline{1}$ (not both) occurring in a fixed position (for a total of $2n$ such cosets).  For brevity, say a coset is ``positive" (resp. ``negative") if it consists of signed permutations having $1$ (resp. $\overline{1}$) in a fixed position.
	\item The minimal length element of such a coset has $2,\hdots,n$ appearing in order in the one-line notation.
	\item The maximal length element of such a coset has $\overline{2},\hdots,\overline{n}$ appearing in order in the one-line notation.
\end{enumerate}

In type $D$, things are similar, but a bit more complicated to state due to the fact that all signed permutations must change an even number of signs:
\begin{enumerate}
	\item As in type $B$, left cosets $W_P \backslash W$ consist of all signed permutations whose one-line notations have either $1$ or $\overline{1}$ (not both) occurring in a fixed position.  Again, we describe these cosets as ``positive" or ``negative" for short.
	\item The minimal length element of such a coset can be described as follows:  If the coset is positive, then the minimal length element has $2,\hdots,n$ appearing in order.  If the coset is negative, then the minimal length element has $2,\hdots,n-1,\overline{n}$ appearing in order.
	\item The description of the maximal length element of such a coset depends on whether $n$ is even or odd.
		\begin{enumerate}
			\item If $n$ is odd, and if the coset is positive, then the maximal length element has $\overline{2},\hdots,\overline{n}$ appearing in order.  If the coset is negative, then the maximal length element has $\overline{2},\hdots,\overline{n-1},n$ appearing in order.
			\item If $n$ is even, and if the coset is positive, then the maximal length element has $\overline{2},\hdots,\overline{n-1},n$ appearing in order.  If the coset is negative, then the maximal length element has $\overline{2},\hdots,\overline{n}$ appearing in order.
		\end{enumerate}
\end{enumerate}

In light of the comments of the introduction, we are interested in Richardson varieties $X_u^v$ where $u$ (resp. $v$) is a maximal (resp. minimal) length coset representative of $W_P \backslash W$.

Given such a $u,v$, say for short that $u$ is ``positive" (resp. ``negative") if it represents a positive (resp. negative) coset, and likewise for $v$.  Then we have the following result on when $u \geq v$ in Bruhat order.  The proof, which we omit, follows easily from the characterizations of the Bruhat order in types $B$ and $D$ given in \cite[\S 3.3,3.5]{Billey-Lakshmibai-00}:

\begin{lemma}\label{lemma:type-bd-bruhat-comparability}
	In type $B$, suppose that $v$ is positive.  If $u$ is negative, then $u \geq v$.  If $u$ is positive, then $u \geq v$ if and only if the $1$ in the one-line notation for $u$ appears to the right of the $1$ in the one-line notation for $v$ (i.e. $u^{-1}(1) \geq v^{-1}(1)$).  If $v$ is negative, then $u \geq v$ if and only if $u$ is also negative, and the $\overline{1}$ in the one-line notation for $u$ appears to the left of the $\overline{1}$ in the one-line notation for $v$ (i.e. $u^{-1}(\overline{1}) \leq v^{-1}(\overline{1})$).  In type $D$, the same holds, except when $v$ is positive, $u$ is negative, and $v^{-1}(1) = u^{-1}(\overline{1}) = n$.  In this case, $u$ and $v$ are not comparable in Bruhat order.
\end{lemma}

For such a pair $u,v$, the Richardson variety $X_u^v$ coincides with an $L$-orbit closure.  To turn this into a Schubert calculus rule, we need to understand the weak order poset of $L$-orbits on $G/B$ and, as mentioned in the introduction, we understand this only conjecturally.  We now state the specifics of Trapa's conjecture regarding the $L$-orbits on $G/B$.

First, recall the parametrization of $K$-orbits given in \cite{Matsuki-Oshima-90}, which is in terms of ``clans".  Generally speaking, a ``$(p,q)$-clan" is a string of $n=p+q$ characters, each a $+$, a $-$, or a natural number.  The natural numbers are paired, meaning that each which appears does so exactly twice, and the number of $+$'s minus the number of $-$'s must be $p-q$.  Such strings are considered equivalent up to permutation of the natural numbers, meaning that it is the \textit{positions} of matching numbers which matters, and not what the numbers actually \textit{are}.  So, for instance, $(1,2,1,2)$, $(2,1,2,1)$, and $(5,7,5,7)$ are all the same clan, while $(1,1,2,2)$ is a different clan.  See \cite{Wyser-12b} for more details if needed.

Here, we will be interested in $(2,2n-1)$-clans in the type $B$ case, and in $(2,2n-2)$-clans in the type $D$ case.  In each case, we will want our clans to satisfy the following symmetry condition.
\begin{definition}
	A clan $\gamma = (c_1,\hdots,c_n)$ is \textit{symmetric} if the clan $\gamma' = (c_n,\hdots,c_1)$ obtained from $\gamma$ by reversing its characters is the same clan.  More explicitly, $\gamma$ is symmetric if and only if
	\begin{enumerate}
		\item If $c_i$ is a sign, then $c_{n+1-i}$ is the same sign.
		\item If $c_i$ is a number, then $c_{n+1-i}$ is also a number, and if $c_i = c_{n+1-j}$, then $c_j = c_{n+1-i}$.
	\end{enumerate}
\end{definition}

By \cite{Matsuki-Oshima-90}, the $K$-orbits are parametrized by symmetric $(2,2n-1)$-clans in the type $B$ case, and by symmetric $(2,2n-2)$-clans in the type $D$ case.  Here is the conjectured combinatorial check on such a clan $\gamma$ which should determine whether the associated $K$-orbit $Q_{\gamma}$ has one or two components.

\begin{conjecture}[P. Trapa]\label{conj:trapa-conjecture}
		For $n$ even or odd, i.e. in both types $B$ and $D$, the $K$-orbit $Q_{\gamma}$ corresponding to the symmetric $(2,n-2)$-clan $\gamma=(c_1,\hdots,c_n)$ is disconnected if and only if for any $i$ with $c_i \in \N$, $c_i \neq c_{n+1-i}$.
\end{conjecture}

Along with this conjecture comes the following description of the weak order on $L \backslash G/B$ relative to the (known) weak order on $K \backslash G/B$:  If, in the weak order graph for $K \backslash G/B$, we have
	\[
		\xymatrixcolsep{1pc} 
		\xymatrixrowsep{1pc}
		\xymatrix
		{
			Q_{\gamma'} \\
			Q_{\gamma} \ar@{->}[u]^\alpha 
		}
	\]
	
	and $Q_{\gamma}$ and $Q_{\gamma'}$ are both disconnected (say splitting as $Q_{\gamma}^1 \cup Q_{\gamma}^2$ and $Q_{\gamma'}^1 \cup Q_{\gamma'}^2$) then in the weak order graph for $L \backslash G/B$, we have
	\[
		\xymatrixcolsep{1pc}
		\xymatrixrowsep{1pc}
		\xymatrix
		{
			Q_{\gamma'}^1 & & Q_{\gamma'}^2 \\
			Q_{\gamma}^1 \ar@{->}[u]^\alpha & & Q_{\gamma}^2 \ar@{->}[u]_\alpha
		}
	\]
	
	As another possibility, if in $K \backslash G/B$ we have
	\[
		\xymatrixcolsep{1pc} 
		\xymatrixrowsep{1pc}
		\xymatrix
		{
			Q_{\gamma'} \\
			Q_{\gamma} \ar@{=>}[u]^\alpha 
		}
	\]
	
	and $Q_{\gamma}$ is disconnected (splitting as $Q_{\gamma}^1 \cup Q_{\gamma}^2$) while $Q_{\gamma'}$ is \textit{connected}, then in the weak order graph for $L \backslash G/B$, we have
	\[
		\xymatrixcolsep{1pc}
		\xymatrixrowsep{1pc}
		
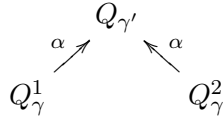
\xymatrix
		{
			& Q_{\gamma'} & \\
			Q_{\gamma}^1 \ar@{->}[ur]^\alpha & & Q_{\gamma}^2 \ar@{->}[ul]_\alpha
		}
	\]
	
	Assuming Conjecture \ref{conj:trapa-conjecture} and examining the known descriptions of the weak order on $K \backslash G/B$ in the cases at hand, it is easy to see that these are the only possibilities.  This thus gives a complete (albeit conjectural) description of $L \backslash G/B$ and its weak order.  For the sake of argument, we assume this combinatorial model of $L \backslash G/B$ to be valid, and deduce from it a Schubert calculus rule describing $c_{w_0u,v}^w$ where $u,v$ are as described in Lemma \ref{lemma:type-bd-bruhat-comparability}.
	
	To deduce the rule, we must describe how to associate to the Richardson variety $X_u^v$ the appropriate $L$-orbit closure.  We do so first in type $B$.  Note that symmetric $(2,2n-1)$-clans are simple to describe.  Indeed, each such clan fits one of the following three descriptions:
	\begin{enumerate}[(a)]
		\item It contains $2$ plus signs, and $2n-1$ minus signs.
		\item It contains $2$ pairs of matching natural numbers, and $2n-3$ minus signs.  The natural numbers may be in any of the patterns $(1,1,2,2)$, $(1,2,1,2)$, or $(1,2,2,1)$.
		\item It contains $1$ pair of matching natural numbers, one plus sign (which must be in the middle position $n+1$), and $2n-2$ minus signs.		
	\end{enumerate}
	
	Having observed this, and recalling the possibilities for $u$ and $v$ stated in Lemma \ref{lemma:type-bd-bruhat-comparability}, we now identify the Richardson variety $X_u^v$ with an $L$-orbit closure.
	\begin{conjecture}\label{conj:type-b-richardson-conj}
	   Let $u,v$ be as described in Lemma \ref{lemma:type-bd-bruhat-comparability}.  We define a clan $\gamma(u,v)$ associated to $u,v$.
		\begin{enumerate}
			\item If $u,v$ are both positive, and $i:=v^{-1}(1) = u^{-1}(1)$, then let $\gamma(u,v)$ be a clan of type (a) above, with the $+$ signs occurring in positions $i$ and $2n+2-i$.
			\item If $u,v$ are both positive, and $i:= v^{-1}(1) < u^{-1}(1) =: j$, then let $\gamma(u,v)$ be a clan of type (b) above, with the natural numbers occurring in positions $i,j,2n+2-j,2n+2-i$, and in the pattern $(1,1,2,2)$.
			\item If $v$ is positive, $u$ is negative, and $i:= v^{-1}(1) \neq u^{-1}(\overline{1}) =: j$, then let $\gamma(u,v)$ be a clan of type (b) above, with the natural numbers occurring in positions $i,j,2n+2-j,2n+2-i$, and in the pattern $(1,2,1,2)$. 
		\suspend{enumerate}
	     Then $X_u^v$ is one of the two $L$-orbit closures associated to the clan $\gamma$, while $X_{w_0v}^{w_0u}$ is the other.
	   
	   On the other hand, 
	   \resume{enumerate}
	   		\item If $v$ is positive, $u$ is negative, and $i:= v^{-1}(1) = u^{-1}(\overline{1}) < n$, let $\gamma(u,v)$ be a clan of type (b), with the $4$ numbers occurring in positions $i,i+1,2n+1-i,2n+2-i$ and in the pattern $(1,2,2,1)$.
	   		\item If $v$ is positive, $u$ is negative, and $v^{-1}(1) = u^{-1}(\overline{1}) = n$, let $\gamma(u,v)$ be a clan of type (c), with the $2$ numbers occurring in positions $n,n+2$.
	   \end{enumerate}
	   Then $X_u^v$ is the single $L$-orbit closure associated to the clan $\gamma$.
	\end{conjecture}
	
	We give examples of the above:
	\begin{enumerate}
		\item $\gamma(\overline{2} 1 \overline{3} \overline{4}, 2134) = (-,+,-,-,-,-,-,+,-)$;
		\item $\gamma(\overline{2} \overline{3} 1 \overline{4}, 2134) = (-,1,1,-,-,-,2,2,-)$;
		\item $\gamma(\overline{2} \overline{3} \overline{1} \overline{4}, 2134) = (-,1,2,-,-,-,1,2,-)$;
		\item $\gamma(\overline{2} \overline{1} \overline{3} \overline{4}, 2134) = (-,1,2,-,-,-,2,1,-)$; 
		\item $\gamma(\overline{2} \overline{3} \overline{4} \overline{1}, 2341) = (-,-,-,1,+,1,-,-,-)$.
	\end{enumerate}
	
We now state the Schubert calculus rule implied by Conjectures \ref{conj:trapa-conjecture} and \ref{conj:type-b-richardson-conj}.  We describe a monoidal action of $W$ on symmetric $(2,2n-1)$-clans, following \cite{Matsuki-Oshima-90}.  Let $\gamma=(c_1,\hdots,c_{2n+1})$ be a symmetric $(2,2n-1)$-clan.
	\begin{definition}\label{def:type-b-action}
		We define the possible ``type $B$ operations" on $\gamma$ as follows:
			\begin{enumerate}[(a)]
				\item Replace $(c_i,c_{i+1})$ by a new pair of matching natural numbers, and $(c_{2n+1-i},c_{2n+2-i})$ by a second new pair of matching natural numbers.
				\item Interchange $(c_i,c_{i+1})$, and also $(c_{2n+1-i},c_{2n+2-i})$.
				\item Interchange $(c_i,c_{i+1})$, \textit{but not} $(c_{2n+1-i},c_{2n+2-i})$.
				\item Interchange $(c_n,c_{n+2})$.
				\item Replace $(c_n,c_{n+2})$ by a new pair of matching natural numbers, and invert the sign in position $c_{n+1}$.
			\end{enumerate}
			
		Then for $i=1,\hdots,n-1$, the action of $s_i$ on $\gamma$ is as follows:
			\begin{enumerate}
				\item If $c_i,c_{i+1}$ are opposite signs, then $s_i \cdot \gamma$ is obtained from $\gamma$ by operation (a).
				\item If $c_i$ is a sign, $c_{i+1}$ is a number, and the mate for $c_{i+1}$ occurs to the right of $c_{i+1}$, then $s_i \cdot \gamma$ is obtained from $\gamma$ by operation (b).
				\item If $c_i$ is a number, $c_{i+1}$ is a sign, and the mate for $c_i$ occurs to the left of $c_i$, then $s_i \cdot \gamma$ is again obtained from $\gamma$ by operation (b).
				\item If $c_i$ and $c_{i+1}$ are unequal natural numbers, with $(c_i,c_{i+1}) = (c_{2n+1-i},c_{2n+2-i})$, then $s_i \cdot \gamma$ is obtained from $\gamma$ by operation (c).
			\suspend{enumerate}
			
		In all other cases, $s_i \cdot \gamma = \gamma$.
		
		The action of $s_n$ on $\gamma$ is as follows:
			\resume{enumerate}
				\item If $c_n$ and $c_{n+2}$ are unequal natural numbers, with the mate for $c_n$ occurring to the left of the mate for $c_{n+2}$, $s_n \cdot \gamma$ is obtained from $\gamma$ by operation (d).
				\item If $(c_n,c_{n+1},c_{n+2}) = (+,-,+)$, $s_n \cdot \gamma$ is obtained from $\gamma$ by operation (e).
				\item If $(c_n,c_{n+1},c_{n+2}) = (-,+,-)$, $s_n \cdot \gamma$ is obtained from $\gamma$ by operation (e). (*)
			\end{enumerate}
		
		In all other cases, $s_n \cdot \gamma = \gamma$.
	\end{definition}
	
	Note the (*) on rule (7); an operation of this sort will introduce a Schubert constant of $2$.
	
	We give examples of rules (1)-(7) above:
	\begin{enumerate}
		\item $s_1 \cdot (+,-,-,-,-,-,+) = (1,1,-,-,-,2,2)$;
		\item $s_1 \cdot (-,1,-,+,-,1,-) = (1,-,-,+,-,-,1)$;
		\item $s_2 \cdot (1,1,-,-,-,2,2) = (1,-,1,-,2,-,2)$;
		\item $s_2 \cdot (-,1,2,-,1,2,-) = (-,1,2,-,2,1,-)$;
		\item $s_3 \cdot (-,1,1,-,2,2,-) = (-,1,2,-,1,2,-)$;
		\item $s_3 \cdot (-,-,+,-,+,-,-) = (-,-,1,+,1,-,-)$.
		\item $s_4 \cdot (-,-,1,-,+,-,1,-,-) = (-,-,1,2,-,2,1,-,-)$.
	\end{enumerate}
	
	Conjectures \ref{conj:trapa-conjecture} and \ref{conj:type-b-richardson-conj} then imply the following rule describing $c_{w_0u,v}^w$:
	\begin{conjecture}\label{conj:type-b-schubert}
		Let $\gamma_0 = (1,2,-,\hdots,-,2,1)$.  Recall rules (1)-(7) given in Definition \ref{def:type-b-action}, which describe the $M(W)$-action.  For $u,v$ as described in the statement of Lemma \ref{lemma:type-bd-bruhat-comparability}, let $\gamma(u,v)$ be the clan defined in the statement of Conjecture \ref{conj:type-b-richardson-conj}.  Then for $w$ of the appropriate length,
		\[ c_{w_0u,v}^w = 
			\begin{cases}
				2 & \text{ if $w \cdot \gamma(u,v) = \gamma_0$ and the computation of the $w$-action involves rule (7),} \\
				1 & \text{ if $w \cdot \gamma(u,v) = \gamma_0$ and the computation of the $w$-action does not involve rule (7),} \\
				0 & \text{ if $w \cdot \gamma(u,v) \neq \gamma_0$.}
			\end{cases}
		\]
	\end{conjecture}
	
	We now describe how to modify things for the type $D$ case.  First, notice that symmetric $(2,2n-2)$-clans are even simpler to describe than symmetric $(2,2n-1)$-clans, as each meets one of the following descriptions:
	
	\begin{enumerate}[(a)]
		\item It contains $2$ plus signs, and $2n-1$ minus signs.
		\item It contains $2$ pairs of matching natural numbers, and $2n-4$ minus signs.  The natural numbers may be in any of the patterns $(1,1,2,2)$, $(1,2,1,2)$, or $(1,2,2,1)$.	
	\end{enumerate}
	
	Then we have the following conjecture regarding type $D$ Richardson varieties $X_u^v$ and $L$-orbit closures:
	\begin{conjecture}\label{conj:type-d-richardson-conj}
		Let $u,v$ be type $D$ Weyl group elements as described in the statement of Lemma \ref{lemma:type-bd-bruhat-comparability}.  Let $\gamma(u,v)$ be the clan associated to $u,v$ as described in the statement of Conjecture \ref{conj:type-b-richardson-conj}, with the exception that we no longer consider case (5).  Then the statement of Conjecture \ref{conj:type-b-richardson-conj} holds also in type $D$.
	\end{conjecture}
	
	\begin{remark}
		Recall that by Lemma \ref{lemma:type-bd-bruhat-comparability}, elements $u,v$ as described in case (5) of Conjecture \ref{conj:type-b-richardson-conj} are not comparable in Bruhat order.  This is why we omit this from consideration in the statement of Conjecture \ref{conj:type-d-richardson-conj}.
	\end{remark}
	
	Assuming Conjecture \ref{conj:type-d-richardson-conj} (which, again, follows easily from Conjecture \ref{conj:trapa-conjecture}), we can state a conjectural rule for Schubert constants $c_{w_0u,v}^w$ in type $D$, but we must define the monoidal $W$-action differently, again following \cite{Matsuki-Oshima-90}.  Let $\gamma=(c_1,\hdots,c_{2n})$ be a symmetric $(2,2n-2)$-clan.
	
	\begin{definition}\label{def:type-d-ops}
		The $M(W)$-action of $s_1,\hdots,s_{n-1}$ in type $D$ is exactly as described in type $B$, cf. Definition \ref{def:type-b-action}.  We define the possible \textbf{type $D$ operations} of $s_n$ on $\gamma$ as follows:
			\begin{enumerate}[(a)]
				\item Interchange $(c_{n-1},c_{n+1})$ and $(c_n,c_{n+2})$.
				\item Replace $(c_{n-1},c_n,c_{n+1},c_{n+2})$ by the pattern $(1,2,1,2)$.
				\item Replace $(c_{n-1},c_n,c_{n+1},c_{n+2})$ by the pattern $(1,2,2,1)$.
			\end{enumerate}
			
			Then the $M(W)$-action on $s_n$ on $\gamma$ is defined as follows:
			\begin{enumerate}
				\item If $(c_{n-1},c_n,c_{n+1},c_{n+2})$ form the pattern $(\pm,1,1,\pm)$, $s_n \cdot \gamma$ is obtained from $\gamma$ by operation (a).
				\item If $(c_{n-1},c_n,c_{n+1},c_{n+2})$ form the pattern $(\pm,1,2,\pm)$, the mate for $c_n$ lies to the left of $c_{n-1}$, and the mate for $c_{n+1}$ lies to the right of $c_{n+2}$, $s_n \cdot \gamma$ is obtained from $\gamma$ by operation (a).
				\item If $(c_{n-1},c_n,c_{n+1},c_{n+2})$ form the pattern $(1,\pm,\pm,2)$, the mate for $c_{n-1}$ lies to the left of $c_{n-1}$, and the mate for $c_{n+2}$ lies to the right of $c_{n+2}$, $s_n \cdot \gamma$ is obtained from $\gamma$ by operation (a).
				\item If $(c_{n-1},c_n,c_{n+1},c_{n+2})$ is equal to $(+,-,-,+)$ or $(-,+,+,-)$, $s_n \cdot \gamma$ is obtained from $\gamma$ by operation (b).
				\item If $(c_{n-1},c_n,c_{n+1},c_{n+2})$ form the pattern $(1,1,2,2)$, $s_n \cdot \gamma$ is obtained from $\gamma$ by operation (c).
			\end{enumerate}
	\end{definition}
	
	We give examples of each of these possibilities:
	\begin{enumerate}
		\item $s_4 \cdot (1,-,-,2,2,-,-,1) = (1,-,2,-,-,2,-,1)$;
		\item $s_4 \cdot (1,-,-,1,2,-,-,2) = (1,-,2,-,-,1,-,2)$;
		\item $s_4 \cdot (1,-,1,-,-,2,-,2) = (1,-,-,2,1,-,-,2)$;
		\item $s_4 \cdot (-,-,+,-,-,+,-,-) = s_4 \cdot (-,-,-,+,+,-,-,-) = (-,-,1,2,1,2,-,-)$;
		\item $s_4 \cdot (-,-,1,1,2,2,-,-) = (-,-,1,2,2,1,-,-)$.
	\end{enumerate}
	
	With the $W$-action so defined, we give the following conjecture in type $D$:
	\begin{conjecture}\label{conj:type-d}
		Let $\gamma_0 = (1,2,-,\hdots,-,2,1)$.  For $u,v$ as described in the statement of Lemma \ref{lemma:type-bd-bruhat-comparability}, let $\gamma(u,v)$ be the clan associated to $u,v$ as described in the statement of Conjecture \ref{conj:type-d-richardson-conj}.  Then for $w$ of the appropriate length,
		\[ c_{w_0u,v}^w = 
			\begin{cases}
				1 & \text{ if $w \cdot \gamma(u,v) = \gamma_0$,} \\
				0 & \text{ otherwise.}
			\end{cases}
		\]
	\end{conjecture}
	
	\begin{example}\label{ex:example-2}
	As an example of a type $B$ Schubert product computed using Conjecture \ref{conj:type-b-schubert}, consider the product $S_{234\overline{1}} \cdot S_{2341}$.  This Schubert product is the class of the Richardson variety $X_u^v$ for $u=\overline{2} \overline{3} \overline{4} 1$ and $v=2341$, to which we associate the clan $\gamma=(-,-,-,+,-,+,-,-,-)$.  We have $l(2341) = 3$ and $l(234 \overline{1}) = 4$, and there are $44$ elements of length $7$ in $W$.  Table 1 shows each of these $44$ elements as words in the simple reflections, the clan obtained from computing the action of each on the clan $\gamma$, and the corresponding structure constant specified by Conjecture \ref{conj:type-b-schubert}.
\end{example}

	\begin{example}\label{ex:example-3}
	As an example of a type $D$ Schubert product computed using Conjecture \ref{conj:type-d}, consider the product $S_{2314} \cdot S_{23\overline{4}\overline{1}}$, the class of the Richardson variety $X_u^v$ with $u = \overline{2} \overline{3} 4 1$ and $v=2314$.  To this Richardson variety, we associate the clan $(-,-,1,1,2,2,-,-)$.  We have $l(2314) = 2$ and $l(23\overline{4}\overline{1}) = 3$, and there are $28$ elements of length $5$ in $W$.  The results of the computation according to Conjecture \ref{conj:type-d} are given in Table 2.
	\end{example}
	
\bibliographystyle{alpha}
\bibliography{../sourceDatabase}

\begin{table}[h]
	\caption{Example \ref{ex:example-2}:  Computing the $B_4$ Schubert product $S_{4321} \cdot S_{432\overline{1}}$}
	\begin{tabular}{|c|c|c|}
		\hline
		Length $7$ Element $w$ & $w \cdot (-,-,-,+,-,+,-,-,-)$ & $c_{u,v}^w$ \\ \hline
		$[1, 2, 1, 4, 3, 2, 1]$ & $(1,-,-,2,-,1,-,-,2)$ & $0$ \\ \hline
        $[3, 2, 1, 4, 3, 2, 1]$ & $(-,1,2,-,-,-,1,2,-)$ & $0$ \\ \hline
        $[1, 3, 2, 4, 3, 2, 1]$ & $(1,-,-,2,-,2,-,-,1)$ & $0$ \\ \hline
        $[2, 3, 2, 4, 3, 2, 1]$ & $(-,1,2,-,-,-,2,1,-)$ & $0$ \\ \hline
        $[2, 1, 3, 4, 3, 2, 1]$ & $(-,1,-,2,-,2,-,1,-)$ & $0$ \\ \hline
        $[1, 2, 3, 4, 3, 2, 1]$ & $(1,-,-,2,-,2,-,-,1)$ & $0$ \\ \hline
        $[1, 3, 2, 1, 4, 3, 2]$ & $(1,-,2,-,-,-,1,-,2)$ & $0$ \\ \hline
        $[2, 3, 2, 1, 4, 3, 2]$ & $(-,1,2,-,-,-,2,1,-)$ & $0$ \\ \hline
        $[4, 3, 2, 1, 4, 3, 2]$ & $(-,1,2,-,-,-,1,2,-)$ & $0$ \\ \hline
        $[2, 1, 3, 2, 4, 3, 2]$ & $(1,2,-,-,-,-,-,1,2)$ & $0$ \\ \hline
        $[1, 2, 3, 2, 4, 3, 2]$ & $(1,-,2,-,-,-,2,-,1)$ & $0$ \\ \hline
        $[1, 2, 1, 3, 4, 3, 2]$ & $(1,-,-,2,-,2,-,-,1)$ & $0$ \\ \hline
        $[2, 1, 3, 2, 1, 4, 3]$ & $(1,2,-,-,-,-,-,1,2)$ & $0$ \\ \hline
        $[1, 2, 3, 2, 1, 4, 3]$ & $(1,-,2,-,-,-,2,-,1)$ & $0$ \\ \hline
        $[1, 4, 3, 2, 1, 4, 3]$ & $(1,-,2,-,-,-,1,-,2)$ & $0$ \\ \hline
        $[2, 4, 3, 2, 1, 4, 3]$ & $(-,1,2,-,-,-,2,1,-)$ & $0$ \\ \hline
        $[3, 4, 3, 2, 1, 4, 3]$ & $(-,1,2,-,-,-,1,2,-)$ & $0$ \\ \hline
        $[1, 2, 1, 3, 2, 4, 3]$ & $(1,2,-,-,-,-,-,2,1)$ & $1$ \\ \hline
        $[2, 1, 4, 3, 2, 4, 3]$ & $(1,2,-,-,-,-,-,1,2)$ & $0$ \\ \hline
        $[1, 2, 4, 3, 2, 4, 3]$ & $(1,-,2,-,-,-,2,-,1)$ & $0$ \\ \hline
        $[3, 2, 4, 3, 2, 4, 3]$ & $(-,1,2,-,-,-,2,1,-)$ & $0$ \\ \hline
        $[1, 3, 4, 3, 2, 4, 3]$ & $(1,-,2,-,-,-,1,-,2)$ & $0$ \\ \hline
        $[2, 3, 4, 3, 2, 4, 3]$ & $(-,1,2,-,-,-,2,1,-)$ & $0$ \\ \hline
        $[1, 2, 1, 3, 2, 1, 4]$ & $(1,-,-,-,+,-,-,-,1)$ & $0$ \\ \hline
        $[2, 1, 4, 3, 2, 1, 4]$ & $(-,1,-,2,-,2,-,1,-)$ & $0$ \\ \hline
        $[1, 2, 4, 3, 2, 1, 4]$ & $(1,-,-,2,-,2,-,-,1)$ & $0$ \\ \hline
        $[3, 2, 4, 3, 2, 1, 4]$ & $(-,1,2,-,-,-,2,1,-)$ & $0$ \\ \hline
        $[1, 3, 4, 3, 2, 1, 4]$ & $(-,-,1,2,-,2,1,-,-)$ & $0$ \\ \hline
        $[2, 3, 4, 3, 2, 1, 4]$ & $(-,1,-,2,-,2,-,1,-)$ & $0$ \\ \hline
        $[1, 2, 1, 4, 3, 2, 4]$ & $(1,-,-,2,-,2,-,-,1)$ & $0$ \\ \hline
        $[3, 2, 1, 4, 3, 2, 4]$ & $(-,1,2,-,-,-,2,1,-)$ & $0$ \\ \hline
        $[1, 3, 2, 4, 3, 2, 4]$ & $(1,-,2,-,-,-,2,-,1)$ & $0$ \\ \hline
        $[2, 3, 2, 4, 3, 2, 4]$ & $(-,1,2,-,-,-,2,1,-)$ & $0$ \\ \hline
        $[2, 1, 3, 4, 3, 2, 4]$ & $(-,1,-,2,-,2,-,1,-)$ & $0$ \\ \hline
        $[1, 2, 3, 4, 3, 2, 4]$ & $(1,-,-,2,-,2,-,-,1)$ & $0$ \\ \hline
        $[1, 3, 2, 1, 4, 3, 4]$ & $(1,-,2,-,-,-,2,-,1)$ & $0$ \\ \hline
        $[2, 3, 2, 1, 4, 3, 4]$ & $(-,1,2,-,-,-,2,1,-)$ & $0$ \\ \hline
        $[4, 3, 2, 1, 4, 3, 4]$ & $(-,1,2,-,-,-,2,1,-)$ & $0$ \\ \hline
        $[2, 1, 3, 2, 4, 3, 4]$ & $(1,2,-,-,-,-,-,2,1)$ & $2$ \\ \hline
        $[1, 2, 3, 2, 4, 3, 4]$ & $(1,-,2,-,-,-,2,-,1)$ & $0$ \\ \hline
        $[1, 4, 3, 2, 4, 3, 4]$ & $(1,-,-,2,-,2,-,-,1)$ & $0$ \\ \hline
        $[2, 4, 3, 2, 4, 3, 4]$ & $(-,1,-,2,-,2,-,1,-)$ & $0$ \\ \hline
        $[3, 4, 3, 2, 4, 3, 4]$ & $(-,1,2,-,-,-,2,1,-)$ & $0$ \\ \hline
        $[1, 2, 1, 3, 4, 3, 4]$ & $(1,-,-,2,-,2,-,-,1)$ & $0$ \\ \hline
	\end{tabular}
\end{table}

\begin{table}[h]
	\caption{Example \ref{ex:example-3}:  Computing the $D_4$ Schubert product $S_{2314} \cdot S_{23\overline{4}\overline{1}}$}
	\begin{tabular}{|c|c|c|}
		\hline
		Length $5$ Element $w$ & $w \cdot (-,-,1,1,2,2,-,-)$ & $c_{u,v}^w$ \\ \hline
		$[2, 1, 3, 2, 1]$ & $(1,-,-,1,2,-,-,2)$ & $0$ \\ \hline
        $[1, 2, 3, 2, 1]$ & $(1,-,-,1,2,-,-,2)$ & $0$ \\ \hline
        $[2, 1, 4, 2, 1]$ & $(1,2,-,-,-,-,1,2)$ & $0$ \\ \hline
        $[1, 2, 4, 2, 1]$ & $(1,-,2,-,-,2,-,1)$ & $0$ \\ \hline
        $[3, 2, 4, 2, 1]$ & $(-,1,2,-,-,2,1,-)$ & $0$ \\ \hline
        $[1, 3, 4, 2, 1]$ & $(1,-,2,-,-,1,-,2)$ & $0$ \\ \hline
        $[2, 3, 4, 2, 1]$ & $(-,1,2,-,-,2,1,-)$ & $0$ \\ \hline
        $[1, 2, 1, 3, 2]$ & $(1,-,-,1,2,-,-,2)$ & $0$ \\ \hline
        $[4, 2, 1, 3, 2]$ & $(1,-,2,-,-,1,-,2)$ & $0$ \\ \hline
        $[1, 2, 1, 4, 2]$ & $(1,2,-,-,-,-,2,1)$ & $1$ \\ \hline
        $[3, 2, 1, 4, 2]$ & $(1,2,-,-,-,-,1,2)$ & $0$ \\ \hline
        $[1, 3, 2, 4, 2]$ & $(1,-,2,-,-,2,-,1)$ & $0$ \\ \hline
        $[2, 3, 2, 4, 2]$ & $(-,1,2,-,-,2,1,-)$ & $0$ \\ \hline
        $[2, 1, 3, 4, 2]$ & $(1,2,-,-,-,-,1,2)$ & $0$ \\ \hline
        $[1, 2, 3, 4, 2]$ & $(1,-,2,-,-,2,-,1)$ & $0$ \\ \hline
        $[1, 4, 2, 1, 3]$ & $(1,-,2,-,-,1,-,2)$ & $0$ \\ \hline
        $[2, 4, 2, 1, 3]$ & $(-,1,2,-,-,2,1,-)$ & $0$ \\ \hline
        $[3, 4, 2, 1, 3]$ & $(-,1,2,-,-,1,2,-)$ & $0$ \\ \hline
        $[2, 1, 4, 2, 3]$ & $(1,2,-,-,-,-,1,2)$ & $0$ \\ \hline
        $[1, 2, 4, 2, 3]$ & $(1,-,2,-,-,2,-,1)$ & $0$ \\ \hline
        $[3, 2, 4, 2, 3]$ & $(-,1,2,-,-,2,1,-)$ & $0$ \\ \hline
        $[1, 3, 4, 2, 3]$ & $(1,-,2,-,-,1,-,2)$ & $0$ \\ \hline
        $[2, 3, 4, 2, 3]$ & $(-,1,2,-,-,2,1,-)$ & $0$ \\ \hline
        $[1, 3, 2, 1, 4]$ & $(1,-,2,-,-,2,-,1)$ & $0$ \\ \hline
        $[2, 3, 2, 1, 4]$ & $(-,1,2,-,-,2,1,-)$ & $0$ \\ \hline
        $[2, 1, 3, 2, 4]$ & $(1,2,-,-,-,-,2,1)$ & $1$ \\ \hline
        $[1, 2, 3, 2, 4]$ & $(1,-,2,-,-,2,-,1)$ & $0$ \\ \hline
        $[1, 2, 1, 3, 4]$ & $(1,-,-,2,2,-,-,1)$ & $0$ \\ \hline
	\end{tabular}
\end{table}

\end{document}